\documentclass{notices}
\usepackage{amsfonts,amssymb,amsmath,amscd,amsthm}
\usepackage{csquotes}
\usepackage{graphicx}
\usepackage{svg}

\numberwithin{equation}{section}
\newtheorem{lemma}{Lemma}[section]
\newtheorem{question}[lemma]{Question}
\newtheorem{theorem}[lemma]{Theorem}
\newtheorem{example}[lemma]{Example}

\newtheorem{definition}[lemma]{Definition}

\title{Hilbert's tenth problem for finitely generated rings}

\author{
  Peter Koymans
  \affil{
    Peter Koymans is a postdoc at Utrecht University. His email address is p.h.koymans@uu.nl.
    }
  \and
  Carlo Pagano
  \affil{
    Carlo Pagano is an assistant professor at Concordia University. His email address is carlein90@gmail.com.
   }
}

\begin{document}
\maketitle

\section{Hilbert's program} 
\label{section1}
In $1900$, Hilbert, encouraged by his friend Minkowski, proposed an extremely influential list of $23$ problems at the ICM in Paris. The tenth problem on the list reads as follows.

\begin{displayquote}
\textbf{Hilbert's tenth problem:} \emph{``Given a Diophantine equation with any number of unknown quantities and with rational integral numerical coefficients: To devise a process according to which it can be determined in a finite number of operations whether the equation is solvable in rational integers.''}
\end{displayquote}

Here, a Diophantine equation refers to an equation of the form 
$$
f(x_1, \ldots, x_n) = 0,
$$
where $f \in \mathbb{Z}[x_1, \ldots, x_n]$ and a solution $(x_1, \ldots, x_n)$ is sought in the set $\mathbb{Z}^n$, as specified at the end of the quoted problem. Long before Hilbert, mathematicians such as Euler, Fermat, and Lagrange, had devised highly ingenious methods to tackle certain special kinds of Diophantine equations. Hilbert's tenth problem informally asks whether the task of solving Diophantine equations could be reduced to an entirely mechanical procedure in complete generality.

\begin{figure}[ht]
\centering
\includegraphics[width=0.3\textwidth]{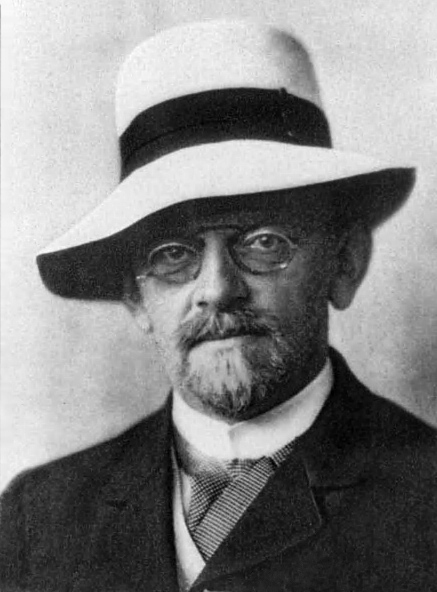}
\caption{David Hilbert}
\end{figure}

Hilbert was well aware of the richness of the theory of Diophantine equations, and his optimism regarding the existence of such a process was a special instance of his general hope of \emph{mechanizing} all of mathematics, also known as \emph{Hilbert's program}. The word ``process'' was intended, in an intuitive sense, to capture what would later be formalized as the notion of an \emph{algorithm}. This soon led to a sequence of dramatic developments that revealed the theoretical limitations of Hilbert's program. 


After these limitations of Hilbert's program came to light, the goal of the mathematical community shifted towards giving a negative answer to Hilbert's tenth problem itself. In $1970$, the combined efforts of Matiyasevich--Robinson--Davis--Putnam culminated in a negative answer for Hilbert's tenth problem. 


Since $1970$, number theorists have sought to establish a negative answer to Hilbert's tenth problem for any infinite number system that is obtained from a finite collection of elements through sum and multiplication. This notion, which is crystallized in the concept of \emph{finitely generated ring} (see Section \ref{section:finitely generated rings}), culminated in the following recent vast generalization of the MRDP result for $\mathbb{Z}$.


\begin{theorem}[{\cite{KP}}]
\label{tKPgen}
Let $R$ be a finitely generated infinite ring. Then Hilbert's tenth problem is undecidable over $R$. 
\end{theorem}

In this article we recount the history of Hilbert's problem, describe its negative solution by MRDP, explain a new approach of Poonen via the theory of elliptic curves, and finally outline the main ideas that go into the proof of Theorem \ref{tKPgen}. In particular, we explain how the proof of Theorem \ref{tKPgen} combines descent, an ancient method in the arithmetic of elliptic curves, with the more modern area of arithmetic combinatorics.


\section{Theoretical limitations} 
\label{section:theoretical limitations}
We recall that Hilbert’s program aimed to mechanize all of mathematics. In particular, Hilbert hoped to give an algorithm for deciding the truth or falsity of any mathematical statement.

In $1931$, G\"odel dealt a first major blow to Hilbert's program, establishing his incompleteness theorem. In his first incompleteness theorem, G\"odel proved that no consistent system of axioms whose theorems can be listed by an algorithm is capable of proving all truths about the arithmetic of natural numbers. 

In $1937$, Turing dealt a second blow to Hilbert's program, showing that there is no algorithm that takes as input any computer program together with an input and decides whether it will halt or run forever on the given input. This is the so-called \emph{halting problem}. In this work, Turing gave solid foundations to the very notion of a computer program, laying the basis for the rise of the digital era. With these foundations in place, we can give a precise formulation of Hilbert's tenth problem.

\begin{displayquote}
\textbf{Precise formulation:} \emph{Is there a Turing machine that takes as input any $f$ in $\mathbb{Z}[x_1, \ldots, x_n]$ and outputs:
\begin{itemize}
    \item YES if $f=0$ has a solution in $\mathbb{Z}^n$,
    \item NO otherwise.
\end{itemize}}
\end{displayquote}

At the heart of both G\"odel's and Turing's proofs is the presence of self-referential statements, an idea that dates back to the liar's paradox \emph{``this statement is false''} from ancient Greece. A closer-in-time predecessor was Cantor's observation that given $\{D_n : n \in \mathbb{Z}\}$, a collection of subsets of $\mathbb{Z}$, the set
$$
\{n \in \mathbb{Z}: n \not \in D_n\}
$$
can not be of the form $D_{n_0}$ for any $n_0$ in $\mathbb{Z}$: if there were such $n_0$, then $n_0$ being in $D_{n_0}$ would force $n_0 \not \in D_{n_0}$, and $n_0$ not being in $D_{n_0}$ would force $n_0 \in D_{n_0}$, a contradiction in both cases. 

Despite these theoretical limitations, Hilbert's program, in its broad aim, is still very much alive. The aim of mechanizing and formalizing mathematical thinking is still a driving force in modern research, especially with the advent of formal verifiers and artificial intelligence. 

\section{Diophantine sets} 
\label{Section: MRDP}
After the work of G\"odel and Turing, the attention of researchers began to shift towards a negative answer to Hilbert's tenth problem. One of the key characters in this line of thinking was Robinson.

\begin{figure}[ht]
\centering
\includegraphics[width=0.3\textwidth]{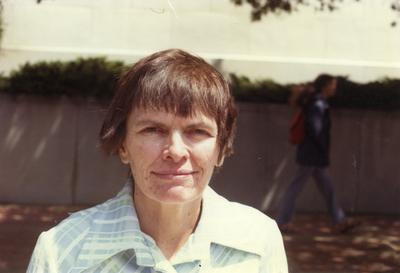}
\caption{Julia Robinson}
\end{figure}

Independently, Davis and later Davis--Putnam began the search for a negative answer. Eventually they joined forces with Robinson, and by $1960$ they were able to give a negative answer to a variant of Hilbert's tenth problem for \emph{exponential} Diophantine equations, namely allowing some of the variables to be exponents in the equation. 

To outline their broad strategy, we introduce a definition. 

\begin{definition}
\label{dDio}
We say that a subset $S \subseteq \mathbb{Z}^n$ is \textup{Diophantine} if there exists some polynomial 
$$
P \in \mathbb{Z}[x_1, \dots, x_n, y_1, \dots, y_m]
$$
such that for all elements $(a_1, \dots, a_n) \in \mathbb{Z}^n$, we have $(a_1, \dots, a_n) \in S$ if and only if there exists $b_1, \dots, b_m \in \mathbb{Z}$ such that 
$$
P(a_1, \dots, a_n, b_1, \dots, b_m) = 0.
$$
\end{definition}

Informally, Diophantine sets are the projection of the vanishing locus of a polynomial equation.

\begin{example}
The set $\mathbb{Z}_{\geq 0}$ of non-negative integers is Diophantine. Indeed, this follows from Lagrange's theorem that every $n \in \mathbb{Z}_{\geq 0}$ can be written as a sum of four integers squares,
\[
n = a^2 + b^2 + c^2 + d^2,
\]
so we can take $P = x_1 - y_1^2 - y_2^2 - y_3^2 - y_4^2$.
\end{example}

\begin{example}
It is typically difficult to prove that a specific subset $S \subseteq \mathbb{Z}^n$ is not Diophantine. However, since $\mathbb{Z}^n$ has uncountably many subsets, while there are only countably many polynomials $P \in \mathbb{Z}[x_1, \dots, x_n, y_1, \dots, y_m]$ (even allowing arbitrary values of $m$), it follows that non-Diophantine subsets must exist.
\end{example}

Likewise, we define a function $f: \mathbb{Z}^n \to \mathbb{Z}$ to be \emph{Diophantine} in case its graph, as a subset of $\mathbb{Z}^{n+1}$, is Diophantine. 

A powerful guiding idea in their approach, familiar from theoretical computer science, is to proceed by reduction to a more classical undecidable problem, such as the halting problem. To achieve this, they had the bold vision to prove that the class of Diophantine functions is precisely the same as the class of functions that can be computed by a Turing machine: their view was that the theory of Diophantine equations is rich enough to simulate essentially any computation. 

This approach can be given the following concrete form. If we could show that every function computable by a Turing machine is in fact Diophantine, then it would be possible to prove that all Diophantine subsets of $\mathbb{Z}$ can be listed as $(D_n)_{n \in \mathbb{Z}}$ such that the set 
\[
U := \{n \in \mathbb{Z} : n \in D_n\} \subseteq \mathbb{Z}
\]
is Diophantine. Let $G(\underline{x}) \in \mathbb{Z}[x_1, \ldots, x_m]$ be a polynomial testifying that $U$ is Diophantine. If there were an algorithm that could check for solutions of the particular family of polynomials
\begin{equation}
\label{Gax}
\{G(a, x_2, \dots, x_m) : a \in \mathbb{Z}\} \subseteq \mathbb{Z}[x_2, \ldots, x_m],
\end{equation}
then it would follow that the complement of $U$ is Diophantine, since its characteristic function would be computable by a Turing machine. This implies that there exists an $n_0 \in \mathbb{Z}$ such that $\mathbb{Z} - U = D_{n_0}$. This allows us to conclude that
$$
D_{n_0} = \mathbb{Z} - U = \{n \in \mathbb{Z}: n \not \in D_n\},
$$
which is impossible by the argument explained at the end of Section \ref{section1}. In their work, Davis--Putnam--Robinson \cite{DPR} showed that it was enough to prove that the function
$$
n \mapsto 2^n
$$
is Diophantine. A slightly more general version of this was known as \emph{hypothesis J.R.} (named after Julia Robinson) and remained open for $20$ years. Despite the immense progress, the entire approach relied on the existence of such exotic Diophantine sets, and there was skepticism in the mathematical community about whether this was possible. 

The situation changed dramatically in $1970$, when the young Matiyasevich \cite{Mat} found precisely one such Diophantine function. 

\begin{figure}[ht]
\centering
\includegraphics[width=0.3\textwidth]{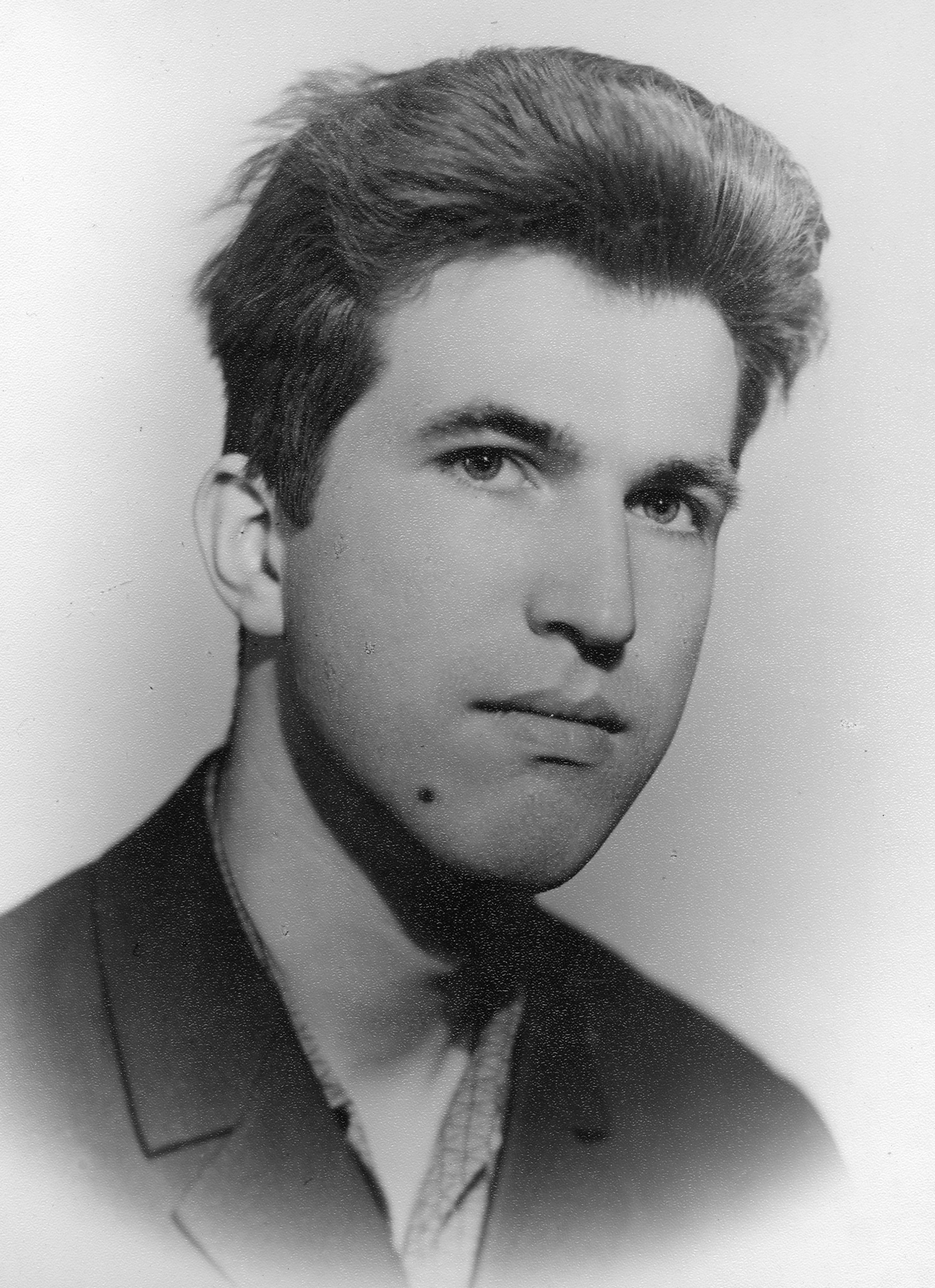}
\caption{Yuri Matiyasevich}
\end{figure}

With an extraordinarily clever argument, he was able to prove that the function
$$
n \mapsto F_n
$$
is Diophantine, where $F_n$ is the $n$-th Fibonacci number. His proof was later streamlined by several authors, and it was identified that the key tool in the argument was the arithmetic of \emph{Pell equations}. A Pell equation is an equation of the form
$$
x^2-dy^2=1,
$$
where $d$ is a positive integer. Subsequent accounts of Matiyasevich's argument exploit the parametric family of Pell equations
\begin{equation}
\label{ePell}
x^2-(a^2-1)y^2=1,
\end{equation}
as $a$ ranges over the integers greater than $1$. Writing, for a positive integer $n$, 
$$
x_n(a) + y_n(a) \cdot \sqrt{a^2 - 1} := (a + \sqrt{a^2-1})^n,
$$
it is not hard to show that all solutions of \eqref{ePell}, with $x, y$ positive integers, are of the form $(x_n(a), y_n(a))$. The key additional property identified by Matiyasevich in gaining access to the exponent $n$ is that 
$$
y_m(a)^2 \mid y_n(a) \Longleftrightarrow y_m(a) \mid n.
$$
Observe that the relation $y_m(a)^2 \mid y_n(a)$ can be reformulated as the equation $y_m(a)^2 \cdot t=y_n(a)$, hence providing Diophantine access to the exponent $n$. Playing this card in an extremely ingenious way, he was able to show that the function
$$
(a, n) \mapsto y_n(a)
$$
is Diophantine. 

From here, we can show that the exponential is Diophantine, and from there conclude, in the way sketched above, that Hilbert's tenth problem is undecidable: there is no algorithm that is capable of deciding, for all Diophantine equations simultaneously, their solvability over $\mathbb{Z}$. All in all, the combined effort of these authors became known as the MRDP theorem. 

\begin{theorem}[``MRDP theorem'']
\label{tMRDP}
There is no algorithm that takes as input an arbitrary $f$ in $\mathbb{Z}[x_1, \ldots, x_n]$ and decides whether $f = 0$ has a solution in $\mathbb{Z}^n$.     
\end{theorem}

By exploiting the polynomial $G$ used to define the set \eqref{Gax}, we see that undecidability holds even if we give a universal bound on the number of variables involved.

Although the MRDP theorem gives a clear theoretical limitation on the existence of a universal algorithm, Hilbert’s tenth problem is still very much alive in spirit. Modern Diophantine geometry is developing a wide range of algorithmic techniques for \emph{special} classes of equations that show up naturally in number theory. For example, decidability is widely believed to hold if one restricts to polynomials $F(x_1, x_2)$ of two variables, i.e., it is believed that there is an algorithm capable of deciding whether an algebraic curve has an integer point.

\section{Finitely generated rings} 
\label{section:finitely generated rings}
After Matiyasevich completed the proof of the MRDP theorem, he raised a natural follow-up question: what happens if, in Hilbert's tenth problem, all instances of $\mathbb{Z}$ are replaced by a finitely generated ring $R$? By definition, a finitely generated ring $R$ is isomorphic to a quotient 
$$
R \cong \mathbb{Z}[x_1, \dots, x_n]/I
$$
for some ideal $I \subseteq \mathbb{Z}[x_1, \dots, x_n]$. This is not only a large and important class of rings, it is a natural class to study for decision problems, since every element of such a ring $R$ can be encoded with a finite amount of data (unlike, for example, the complex numbers).

\begin{example}
The usual integers $\mathbb{Z}$, the Gaussian integers $\mathbb{Z}[\sqrt{-1}] := \{a + b \sqrt{-1} : a, b \in \mathbb{Z}\}$, and the polynomial rings $\mathbb{Z}[x]$ and $\mathbb{Z}[x, y]$ are all finitely generated rings.
\end{example}

\begin{example}
\label{e1}
Examples of rings that are not finitely generated are the rational numbers $\mathbb{Q}$, the real numbers $\mathbb{R}$ or the complex numbers $\mathbb{C}$.
\end{example}

Note that a polynomial equation $f(x_1, \dots, x_n) = 0$ over a finite ring $R$ can easily be solved by evaluating the polynomial at every $(a_1, \ldots, a_n) \in R^n$. We henceforth restrict our attention to the case of finitely generated rings of infinite cardinality.


It is reasonable to expect that the answer to Hilbert's tenth problem should still be negative for such rings. It turns out that the reduction strategy again rears its head, where we now try to reduce Hilbert’s tenth problem over $R$ to Hilbert's tenth problem over $\mathbb{Z}$. Indeed, by encoding a Diophantine equation with coefficients in one ring inside a Diophantine equation having coefficients in another ring, we can hope to transfer undecidability from a known case to many new settings. This leads to the following generalization of Definition \ref{dDio} in which the ring $\mathbb{Z}$ is replaced with the ring $R$.

\begin{definition}
We say that a subset $S \subseteq R^n$ is $R$-Diophantine if there exists some polynomial 
$$
P \in R[x_1, \dots, x_n, y_1, \dots, y_m]
$$
such that $(a_1, \dots, a_n) \in S$ if and only if there exists $b_1, \dots, b_m \in R$ such that 
$$
P(a_1, \dots, a_n, b_1, \dots, b_m) = 0.
$$
\end{definition}

In other words, an $R$-Diophantine set is the projection onto the first $n$ coordinates of the solution set of a polynomial equation over $R$. This notion is central to reduction arguments for Hilbert’s tenth problem. Indeed, if the set of integers $\mathbb{Z}$ can be realized as an $R$-Diophantine subset of $R$, then every Diophantine equation over $\mathbb{Z}$ can be encoded as a Diophantine equation over $R$. It follows that Hilbert’s tenth problem over $R$ is undecidable.

So how can we prove that $\mathbb{Z}$ is Diophantine over $R$? This turned out to be a difficult task, and after the proof of MRDP in 1970, a complicated web of reduction steps started to appear in the literature. Eisentr\"ager \cite{Eisentrager} proved in her PhD thesis that, in order to settle Hilbert's tenth problem for all finitely generated rings, it suffices to prove that $\mathbb{Z}$ is Diophantine over a much smaller and more tractable class of rings, as we now explain.

\begin{definition}
An integral domain $R$ is called an \emph{order} if its underlying additive group is isomorphic to $\mathbb{Z}^n$ for some $n \geq 1$.
\end{definition}

\begin{example}
Every order is a finitely generated ring. The converse, however, is false. The polynomial ring $\mathbb{Z}[x]$ is finitely generated as a ring, but it is not an order, because its underlying additive group is isomorphic to a direct sum of countably many copies of $\mathbb{Z}$. Examples of orders include $\mathbb{Z}$ itself, the ring of Gaussian integers $\mathbb{Z[}\sqrt{-1}]$, and the finite index subring of the Gaussian integers defined by
$$
\left\{a + 2b \sqrt{-1} : a, b \in \mathbb{Z}\right\}.
$$
\end{example}

In particular, Eisentr\"ager’s work shows that it suffices to restrict attention to \emph{maximal orders} $R$; that is, orders for which there is no strictly larger order $S$ with
$$
R \subsetneq S \quad \text{ and } \quad [S : R] < \infty.
$$
Maximal orders enjoy especially nice algebraic properties and play a central role in number theory. In many regards, they are the arithmetic analogues of non-singular curves in algebraic geometry. Maximal orders are often called \emph{rings of integers}, reflecting their close resemblance to the usual integers $\mathbb{Z}$.

Even earlier, Denef \cite{Denef}, building on ideas of Matiyasevich and exploiting the arithmetic of Pell equations, proved the undecidability of Hilbert’s tenth problem for many such maximal orders. However, Denef's approach applies only when the maximal order $R$ is \emph{totally real}, meaning that every embedding of $R$ into $\mathbb{C}$ actually has image contained in $\mathbb{R}$. Although these results represented significant progress, a large class of maximal orders lay beyond the reach of Pell-type methods. For these remaining cases, a fundamentally new idea was needed.

\section{Elliptic curves} 
\label{section:elliptic curves}
An elliptic curve $E$ over the rational numbers $\mathbb{Q}$ is, by definition, an equation of the form
\begin{equation}
\label{ellipticcurve}
E: y^2 = x^3 + ax + b
\end{equation}
where $a, b \in \mathbb{Q}$ are fixed and the cubic on the right-hand side has no repeated roots.

\begin{figure}[ht]
\centering
\includegraphics[width=0.45\textwidth]{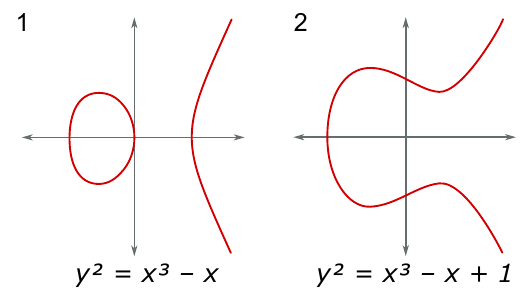}
\caption{Two elliptic curves depicted over the real numbers}
\end{figure}

Elliptic curves have taken a prominent role inside number theory. Most famously, they were a key ingredient in Wiles’s proof of Fermat’s Last Theorem. Despite decades of intensive study, many aspects of elliptic curves remain deeply mysterious.

Already in ancient Greece, it was understood that solutions to \eqref{ellipticcurve} could be used to create new solutions, which in modern terms leads to the fact that an elliptic curve has a natural group structure. To define it, we first pass to the projective closure of the curve, which amounts to adjoining a single additional point, called the point at infinity and denoted by $O$. This point serves as the identity element of the group.

\begin{figure}[ht]
\centering
\includegraphics[width=0.3\textwidth]{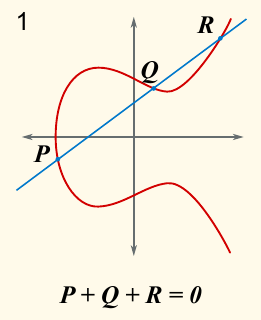}
\caption{The addition law}
\label{f3}
\end{figure}

The group law admits a simple geometric description. If three points $P$, $Q$ and $R$ lie on the elliptic curve and on a straight line, then we declare 
$$
P + Q + R = O. 
$$
For an illustration of the addition law, see Figure \ref{f3}.

Given an elliptic curve $E$, we define its set of \emph{rational points} by
$$
E(\mathbb{Q}) := \{(x, y) \in \mathbb{Q} \times \mathbb{Q} : y^2 = x^3 + ax + b\} \cup \{O\}.
$$
More generally, we can define $E(K)$ for any field $K$.

In the early 1900s, Mordell proved the following fundamental result:

\begin{theorem}
\label{tMordell}
We have 
$$
E(\mathbb{Q}) \cong \mathbb{Z}^r \times \textup{(finite group)}.
$$
The integer $r$ is called the \emph{rank} of the elliptic curve.
\end{theorem}

This makes $E$ a natural candidate to replace the role of Pell's equation in Denef's argument. Note that the analogy with the Pell equation is particularly close when the rank of $E$ equals $1$; in that case, the solution set to both equations have a natural group structure with the same rank. However, elliptic curves come with one rather serious drawback. Despite the simple form that the abelian group $E(\mathbb{Q})$ has, there is currently no known (provably correct) algorithm to compute the rank $r$. And correspondingly, the proof of Theorem \ref{tMordell} is non-constructive.

So far for the bad news. The first piece of good news is that Poonen \cite{Poonen} was able to replace the role of Pell's equation in Denef's arguments by an elliptic curve. This led to the following theorem.

\begin{theorem}[{\cite{Poonen}}]
\label{tElliptic}
Let $R$ be a maximal order. Let $K$ be the field of fractions of $R$. Suppose that there exists an elliptic curve $E: y^2 = x^3 + ax + b$ with $a, b \in \mathbb{Q}$ such that
\begin{equation}
\label{ePoonen}
\mathrm{rank} \, E(\mathbb{Q}) = \mathrm{rank} \, E(K) = 1. 
\end{equation}
Then $\mathbb{Z}$ is $R$-Diophantine. 
\end{theorem}

At the time of Poonen’s work, even finding elliptic curves with $\mathrm{rank} \, E(K) = 1$ was an open problem. Nevertheless, Theorem \ref{tElliptic} marked a major conceptual breakthrough. It provided a clear roadmap for showing that $\mathbb{Z}$ is $R$-Diophantine, and thus for proving that Hilbert’s tenth problem has a negative answer for all finitely generated rings $R$ of infinite cardinality.

As mentioned previously, there is no known algorithm to compute the rank. However, the situation changes drastically if we are willing to assume the Birch--Swinnerton-Dyer (BSD) conjecture. This conjecture is one of the Millennium Problems and is currently wide open. One amazing consequence of the BSD conjecture is that the rank of an elliptic curve should equal the order of vanishing of a certain complex analytic function; in particular, this gives a concrete way to compute ranks. Under this assumption, Mazur and Rubin \cite{MR}, building on an extension of Poonen's result due to Shlapentokh \cite{Shl}, proved the following:

\begin{theorem}[{\cite{MR}}]
Assume the BSD conjecture. Let $R$ be a maximal order. Then $\mathbb{Z}$ is $R$-Diophantine. In particular, Hilbert's tenth problem is undecidable over $R$.
\end{theorem}

Combining this with Eisentr\"ager’s result, we deduce that, \emph{assuming BSD}, Hilbert’s tenth problem is undecidable for every finitely generated ring $R$ of infinite cardinality.

\section{Recent developments} \label{section:recent developments}
Following the work of Mazur and Rubin, mathematicians began searching for an unconditional resolution of Hilbert’s tenth problem for finitely generated rings. To start, let us put aside the additional complications introduced by equation \eqref{ePoonen} and focus on the core challenge: proving the existence of a rank $1$ elliptic curve over $K$. 

A key tool in this quest, already employed by Mazur and Rubin, is \emph{descent}. Descent is a classical method in number theory, with roots dating back to ancient Greece. Fermat himself mastered descent, famously using it to settle the case $n = 4$ of Fermat's Last Theorem. Today, descent remains a powerful method, although its modern incarnations are significantly more abstract. 

Here is an overview of how descent works. Suppose that we want to prove that a certain equation $f(x, y) = 0$ has no solutions in positive integers. What we do is show that if someone gives us a solution $(x_1, y_1)$, we can use the given solution to create a new ``smaller'' solution $(x_2, y_2)$ satisfying
\[
x_2 < x_1. 
\]
Note the strict inequality. Repeating this process, we see that even a single solution can be used to create an infinite list of solutions with strictly decreasing positive $x$-coordinates. This contradiction shows that there are no solutions! How, you might ask, do we use the given solution to create a smaller new solution? The answer varies from problem to problem, ranging from a relatively easy argument used by the Pythagorean school to prove that $\sqrt{2}$ is irrational to much deeper methods used to prove Theorem \ref{tMordell}.

Let $R$ be an order, let $K$ be its field of fractions, and let $E$ be an elliptic curve defined by an equation with coefficients in $K$. When applied to the elliptic curve $E$, the descent machinery provides a sequence of Selmer groups $S(E, K, n)$, one for each $n \geq 2$. The group $S(E, K, n)$ has the following properties:
\begin{itemize}
\item
$S(E, K, n)$ is effectively computable.
\item
$S(E, K, n)$ is more-or-less a product of cyclic groups of order $n$. We denote its rank (number of generators) by $s(E, K, n)$.
\item The quotient group $E(K)/nE(K)$ sits as a subgroup of $S(E, K, n)$. In particular, for all $n \geq 2$,
\[
\text{rank of $E(K)$} \leq s(E, K, n),
\]
where we emphasize that for any given~$n$, the \emph{Selmer rank} $s(E, K, n)$ is effectively computable.
\end{itemize}
Assuming the BSD conjecture, the minimum of $s(E, K, p^k)$ over all $k$ equals the rank of $E(K)$. Moreover, under mild hypotheses, $s(E, K, p^k)$ has the same parity as the rank of $E(K)$. In practice, even $s(E, K, 2)$ frequently equals the rank of $E$.

In several notable instances, it is possible to get a handle on the statistical behavior of $s(E, K, 2)$ as $E$ varies through a family of elliptic curves. In fact, Mazur--Rubin \cite{MR} show that for every maximal order $R$, there exists an elliptic curve $E$ over $R$ such that $s(E, K, 2) = 1$. Since $s(E, K, 2)$ is an upper bound for the rank of $E$ and, assuming BSD, has the same parity, it follows that the rank of $E$ must also be $1$. This observation lies at the heart of their argument. In order to make this argument unconditional, we must simultaneously guarantee two separate ingredients:
\begin{enumerate}
\item[$(i)$] find a family $\mathcal{F}$ of elliptic curves for which every $E \in \mathcal{F}$ has rank at least $1$,
\item[$(ii)$] the existence of an elliptic curve $E \in \mathcal{F}$ with $s(E, K, 2) \leq 1$.
\end{enumerate}
Condition $(i)$ provides a uniform lower bound on the rank, while $(ii)$ furnishes a corresponding upper bound. Together, they imply the existence of an elliptic curve of rank exactly $1$.

Individually, each of the conditions $(i)$ and $(ii)$ is not too difficult to achieve. For example, for $(i)$, we choose the coefficients $a$ and $b$ in a clever way that guarantees a rational point of infinite order. Condition $(ii)$, on the other hand, is precisely what Mazur and Rubin established for certain carefully chosen families. What had been missing, until very recently, was a method for controlling $s(E, K, 2)$ in families $\mathcal{F}$ for which every member $E \in \mathcal{F}$ is known a priori to have rank at least $1$.

The first work to successfully combine $(i)$ and $(ii)$ is due to the authors \cite{KP}. Their starting point is the construction of a suitable family of elliptic curves. Fix elements $a_1, a_2, a_3 \in R$, and consider the elliptic curve
$$
\widetilde{E} : y^2 = (x - a_1) (x - a_2) (x - a_3) =: f(x),
$$
whose defining cubic\footnote{This special form plays a crucial role in the descent arguments.} $f(x)$ splits completely into linear factors.

From this, they define a two-parameter family $\mathcal{F} = (E_{m, n})_{m, n}$ by
$$
y^2 = \bigl(x - a_1 f(m/n)\bigr) \bigl(x - a_2 f(m/n)\bigr) \bigl(x - a_3 f(m/n)\bigr).
$$
Each curve $E_{m, n}$ in this family comes equipped with a tautological rational point, given by
$$
x = f(m/n) \cdot m/n, \qquad y = f(m/n),
$$
so condition $(i)$ is automatically satisfied. The remaining challenge is to establish condition $(ii)$: how can we control the $2$-Selmer rank $s(E_{m, n}, K, 2)$ in a family where a rational point is built in from the outset?

It is classically known that the Selmer ranks $s(E, K, 2)$ in the family $\mathcal{F}$ can be computed explicitly in terms of the prime divisors of the discriminant of the cubic polynomial defining the elliptic curve, namely
$$
(x - a_1 f(m/n)) (x - a_2 f(m/n)) (x - a_3 f(m/n)). 
$$
Because this cubic splits into three linear factors and the coefficients $a_1, a_2, a_3$ are fixed, the main contribution to this discriminant is the quantity
\begin{equation}
\label{nma}
n (m - a_1n) (m - a_2n) (m - a_3n).
\end{equation}
The family $\mathcal{F}$ was engineered so that the expression \eqref{nma} is a product of four linear forms. As a consequence, all four factors can be made simultaneously prime by the Green–Tao theorem, recently generalized to maximal orders by Kai \cite{Kai}. Using this, the authors prove:

\begin{theorem}[{\cite{KP}}]
\label{tKP}
Let $R$ be a maximal order. Then $\mathbb{Z}$ is $R$-Diophantine. In particular, Hilbert's tenth problem is undecidable over $R$. 
\end{theorem}

As an immediate consequence, Hilbert’s tenth problem is undecidable over every finitely generated ring of infinite cardinality.

Some months after the release of \cite{KP}, two other groups of authors \cites{ABHS, Zywina} independently resolved Hilbert's tenth problem for finitely generated rings. While the specific families and descent arguments in these papers differ substantially in their technical details, all three approaches share a common underlying strategy by combining an arithmetic combinatorics result with descent methods. It is worth mentioning that the works \cites{ABHS, Zywina} require only three term progressions in the primes, while \cite{KP} needs a more advanced result.

Despite the recent progress, there is one major open problem in the field that the authors would like to advertise. As already noted in Example \ref{e1}, the field $\mathbb{Q}$ is not finitely generated, and so Theorem \ref{tKP} does not apply. Indeed, Hilbert’s tenth problem over $\mathbb{Q}$ remains wide open:

\begin{question}
\label{qHilbertQ}
Determine whether there exists an algorithm that on input a polynomial $f \in \mathbb{Q}[x_1, \dots, x_n]$ outputs:
\begin{itemize}
    \item \emph{YES} if $f$ has a zero $(a_1, \dots, a_n) \in \mathbb{Q}^n$, and
    \item \emph{NO} otherwise.
\end{itemize} 
\end{question}

It is generally believed that the answer to Question \ref{qHilbertQ} should be no.

In stark contrast to the finitely generated setting, we currently do not know sufficient conditions to answer Question \ref{qHilbertQ} in the negative, nor is there even a clear strategy. Identifying a plausible route toward answering  Question \ref{qHilbertQ} appears to be a formidable challenge!

\bibliography{ExampleRefs}
\end{document}